\documentclass[12pt]{elsarticle}
\usepackage{graphics}
\usepackage{graphicx}
\usepackage{subfigure}

\usepackage{latexsym}
\usepackage{amsmath}
\usepackage{graphicx}
\usepackage{times}
\usepackage{graphicx,multicol,enumerate,subfigure,amsmath}
\usepackage{float}
\usepackage{fullpage}
\usepackage{setspace}
\usepackage[english]{babel}
\usepackage[version=3]{mhchem}
\usepackage{amsmath, amssymb, amsthm}
\usepackage{verbatim, latexsym}
\usepackage{colortbl}

\theoremstyle{definition}
\newtheorem{exmp}{Example}[section]

\usepackage{soul}
\setstcolor{red}

\biboptions{sort&compress}
\journal{}

\onehalfspace

\begin{document}

%\begin{frontmatter}

\title{ Global-Local Nonlinear Model Reduction for Flows in Heterogeneous Porous Media\\
 {\it\small Dedicated to Mary Wheeler on the occasion of her 75-th
  birthday anniversary}\\
}
\author{\textbf{Manal Alotaibi}$^1$}
\author{\textbf{Victor M. Calo}$^{2,3}$, \textbf{Yalchin
    Efendiev}$^{1,2}$\corref{cor1}}
\cortext[cor1]{Email address : efendiev@math.tamu.edu; Telephone: +1-979-458-0836; Fax: +1-979-862-4190}

\author{\textbf{Juan Galvis}$^4$ }
\author{\textbf{Mehdi Ghommem}$^2$ }

\address{$^{1}$
Department of Mathematics \& Institute for Scientific Computation (ISC) \\
  Texas A\&M University \\
  College Station, Texas, USA}

\address{$^{2}$
  Center for Numerical Porous Media (NumPor) \\
  King Abdullah University of Science and Technology (KAUST) \\
  Thuwal 23955-6900, Kingdom of Saudi Arabia}

\address{$^{3}$
Applied Mathematics \& Computational Science and Earth
  Sciences \& Engineering \\
  King Abdullah University of Science and Technology (KAUST) \\
  Thuwal 23955-6900, Kingdom of Saudi Arabia}

\address{$^{4}$
  Departamento de Matem\'aticas, Universidad Nacional de Colombia,\\
  Carrera 45 No 26-85 - Edificio Uriel Gutierr\'ez, Bogot\'a D.C. -
  Colombia}

\begin{abstract}

  In this paper, we combine discrete empirical interpolation
  techniques, global mode decomposition methods, and local multiscale
  methods, such as the Generalized Multiscale Finite Element Method
  (GMsFEM), to reduce the computational complexity associated with
  nonlinear flows in highly-heterogeneous porous media. To solve the
  nonlinear governing equations, we employ the GMsFEM to represent the
  solution on a coarse grid with multiscale basis functions and apply
  proper orthogonal decomposition on a coarse grid.  Computing the
  GMsFEM solution involves calculating the residual and the Jacobian
  on the fine grid. As such, we use local and global empirical
  interpolation concepts to circumvent performing these computations
  on the fine grid. The resulting reduced-order approach enables a
  significant reduction in the flow problem size while accurately
  capturing the behavior of fully-resolved solutions.  We consider
  several numerical examples of nonlinear multiscale partial
  differential equations that are numerically integrated using fully-
  implicit time marching schemes to demonstrate the capability of the
  proposed model reduction approach to speed up simulations of
  nonlinear flows in high-contrast porous media.

\begin{keyword}
  Generalized multiscale finite element method, Nonlinear PDEs,
  heterogeneous porous media, discrete empirical interpolation, proper
  orthogonal decomposition.
\end{keyword}

\end{abstract}

\maketitle
%\tableofcontents

%====================

\section{Introduction}
Nonlinear partial differential equations (PDEs), with multiple scales
and/or high contrast in media properties, represent a class of
problems with many relevant engineering and scientific applications in
porous media. Solving these equations using iterative methods, such as
Newton iterations, requires updating the numerical solution of a large
system of equations at each iteration using the previous iterate
results. Due to the extensive computational requirements
resulting from the disparity of scales and nonlinearity, computing
fine-grid solutions becomes prohibitively expensive. Moreover, these
type of problems often involve coefficients that exhibit
high-contrast and heterogeneous distributions. For
example, when modeling subsurface flows the underlying permeability
field is often represented as a high- contrast coefficient in the
pressure equation. This complicates the simulation of a
large number of configurations for design purposes. As such, we
develop simplified reduced-order models to speed up simulations of
nonlinear flows in porous media within a prescribed accuracy. In this
paper, we present a coupled approach for solving nonlinear PDEs that
combines local and global model-reduction techniques and
the Discrete Empirical Interpolation Method (DEIM).

Proper Orthogonal Decomposition (POD) is one of the best
known global model-reduction methods. The main purpose of this
technique is to reduce the dimension of the dynamical system by
projecting the high-dimensional system into a lower-dimensional
manifold using a set of orthonormal basis functions (POD modes)
constructed from a sequence of snapshots~\cite{Akhtar2009B, Wang2011A,
  Wang2012A, Akhtar2012A, GhommemPCFD2012}. In addition to order
reduction, this technique constitutes a powerful mode decomposition
technique for extracting the most energetic structures from a linear
or nonlinear dynamical process~\cite{Lumley1967, Sirovich1987A,
  Deane1991, Berkooz1993, Holmes1996, Akhtar2009B, Wang2011A,
  Wang2012A, Akhtar2012A, Hay2010A, Hay2011}.

Coarse-grid computational models are often preferred
  because of the computational cost of solving the systems arising in
  the approximation of the nonlinear flow equation on the fine grid.
Some accurate reduced-order methods have been introduced and used in
various applications, such as Galerkin multiscale finite elements
(e.g.,~\cite{Arbogast_two_scale_04, Chu_Hou_MathComp_10,
  egw11,eh09,ehg04, GhommemJCP2013}), mixed multiscale finite element
methods (e.g.,~\cite{aarnes04, ae07, Arbogast_Boyd_06,Iliev_MMS_11}),
the multiscale finite volume method (see, e.g.,~\cite{jennylt03}), and
mortar multiscale methods (see e.g.,~\cite{Arbogast_PWY_07,
  Wheeler_mortar_MS_12}), where Mary Wheeler and her collaborators
have made significant contributions. The main idea of the these
methods is to construct coarse basis functions that approximate the
solution on a coarse grid. Multiscale methods can be considered as
local model reduction techniques that approximate the solution on a
coarse grid for arbitrary coarse- level inputs. In this paper we apply
the enriched coarse space construction from the Generalized Multiscale
Finite Element Method (GMsFEM) as an effective tool for local model
reduction~\cite{EGG_MultiscaleMOR, egw11, ge10, ge10part2}.

The combining of the aforementioned local and global model-reduction
schemes has been used for linear problem~\cite{GhommemJCP2013,
  EGG_MultiscaleMOR}. A significant reduction in the computational
complexity when solving linear parabolic PDE in~\cite{GhommemJCP2013}
has been achieved by combining the concepts of (GMsFEM) and (POD)
and/or Dynamic Mode Decomposition (DMD).  In~\cite{EGG_MultiscaleMOR},
balanced truncation is used to perform global model reduction and is
efficiently combined with the local model reduction tools introduced
in~\cite{egw11}. More recently, local and global multiscale methods
are combined to derive reduced-order models for nonlinear flows in
high-contrast porous media. In~\cite{Ghommem_nonlinPDEs_JCP2013}, the
proposed multiscale empirical interpolation method for solving
nonlinear multiscale PDEs uses GMsFEM to represent the coarse-scale
solution. To avoid performing fine-grid computations, the
discrete empirical interpolation method introduced in~\cite{cs10} was
used to approximate the nonlinear functions at selected points in each
coarse region and then a multiscale proper orthogonal decomposition
technique is used to find an appropriate interpolation
vector. Although, the numerical results presented
in~\cite{Ghommem_nonlinPDEs_JCP2013} proved the applicability of the
presented method, the reduction is limited by the full cost of the
evaluation of the projected nonlinear function.  When
dealing with reduced-order models of nonlinear systems obtained by
projecting the governing equations onto a subspace spanned by the POD
modes, the evaluation of the projected nonlinear term is costly since
it depends on the full dimension of the original system.

In this paper, our main contribution is to circumvent this issue by
employing DEIM to approximate the nonlinear functions locally (at
selected points in each coarse region) at the offline stage and
globally (at selected points in the domain) at the online stage. For
this reason, we refer to our method as global-local
nonlinear approach. The numerical results presented in
this paper show that the proposed method enables
significant reduction in the computational cost associated with
constructing projection-based reduced-order models. In addition to the
model reduction, the proposed approach allows us to improve the
reduced-order solutions in different ways. For instance, increasing
the number of local and global points used at the offline and online
stages, respectively, leads to better approximation (see
Example~\ref{ex2}). Also, using several offline parameter inputs
(Example~\ref{ex3}) improves the reduced-order solutions.

The remainder of the paper is organized as follows. In
Section~\ref{sec:Preliminaries} we introduce and describe the model
problem, the discrete empirical interpolation method, and generalized
multiscale finite element method. The presented global-local DEIM
approach is then discussed in Section~\ref{sec:globallocal}. Numerical
results are presented in Section~\ref{sec:num} and conclusions
in~\ref{sec:conc}.

\section{Preliminaries}
\label{sec:Preliminaries}
\subsection{Model problem}
\label{modproblem}
%====================================
We consider a time-dependent nonlinear flow governed by the following
parabolic partial differential equation
\begin{equation} \label{parabolic}
\frac{\partial u}{\partial t} - \nabla \cdot \left( \kappa(x;u,\mu) \nabla u \right) = h(x)  \quad \text{in} ~ \Omega,
\end{equation}
{with some boundary conditions. The variable } $u= u(t,x; \mu)$ denotes the pressure, $\Omega$ is a bounded
domain, $h$ is a forcing term, and in our case the permeability field
represented by $\kappa(x;u,\mu)$ is a nonlinear function. Here,
$\frac{\partial }{\partial t}$ is the time derivative and $\mu$
represents a parameter.
%, $\nabla=(\frac{\partial }{\partial x},\frac{\partial }{\partial y})$, and $\nabla \cdot \nabla = \Delta = \frac{\partial^2 }{\partial x^2}+\frac{\partial^2 }{\partial y^2}$.

\subsection{Discrete empirical interpolation method (DEIM)}
\label{SEC:DEIM}
We approximate with the Discrete Empirical Interpolation Method
(DEIM)~\cite{cs10} local and global nonlinear functions. DEIM is based
on approximating a nonlinear function by means of an interpolatory
projection of a few selected snapshots of the function. The idea is to
represent a function over the domain while using empirical snapshots
and information in some locations (or components).

We briefly review DEIM as presented in~\cite{cs10}.  Let $ {f}(\tau)
\in \mathbb{R}^{n}$ denotes a nonlinear function where $\tau \in
\mathbb{R}^{n_s}$.  Here, in general, $n_s$ can be different from
$n$. In a reduced-order modeling, $\tau$ has a reduced representation
\begin{align*}
\tau = \sum_{i=1}^{l} \alpha_i \zeta _i
\end{align*}
where $l \ll n_s$. This leads us to look for an approximation of
$f(\tau)$ at a reduced cost. To perform a reduced order approximation
of $f(\tau)$, we first define a reduced dimensional space for
$f(\tau)$. That is, we would like to find $m$ basis vectors (where $m$
is much smaller than $n$), $\psi_1$,..., $\psi_m$, such that we can
write
\begin{equation}\label{fapp}
  {f}(\tau) \approx  \Psi  {d}(\tau),
\end{equation}
where $ \Psi=( \psi_{_1},\cdots, \psi_{_m}) \in
\mathbb{R}^{n\times m}$.

The goal of DEIM is to find $d(\tau)$ using only a few rows
of~\eqref{fapp}.  In general, one can define $d(\tau)$'s using $m$
rows of~\eqref{fapp} and invert a reduced system to compute
$d(\tau)$. This can be formalized using the matrix $\mbox{P}$
$$
\mbox{P} = [ {e}_{\wp_1},\cdots, {e}_{\wp_m}] \in \mathbb{R}^{n\times
  m},
$$
where $ {e}_{\wp_i}=[0,\cdots,0,1,0,\cdots,0]^T\in \mathbb{R}^{n}$ is
the $\wp_i^{\mbox{th}}$ column of the identity matrix $\mbox{I}_n \in
\mathbb{R}^{n\times n}$ for $i=1,\cdots,m$. Multiplying
Equation~\eqref{fapp} by $ \mbox{P}^T $ and assuming that the matrix
$\mbox{P}^T \Psi$ is nonsingular, we obtain
\begin{equation}\label{fapp2}
  {f}(\tau) \approx \tilde{f}(\tau) = \Psi  d(\tau)= \Psi (\mbox{P}^T
  \Psi)^{-1}\mbox{P}^T {f}(\tau).
\end{equation}

To summarize, approximating the nonlinear function $ {f}(\tau)$, as
given by Equation~\eqref{fapp2}, requires the following:
\begin{itemize}
\item Computing the projection basis $ \Psi=( \psi_{_1},\cdots,
  \psi_{_m})$;
\item Identifying the indices $\{\wp_1,\cdots,\wp_m\}.$
\end{itemize}

To determine the projection basis $ \Psi=( \psi_{_1},\cdots,
\psi_{_m})$, we collect function evaluations in an $n \times n_s$
matrix $\mbox{F}=[ {f}(\tau_1),\cdots, {f}(\tau_{n_s})]$ and employ
POD to select the most energetic modes.  This selection uses the
eigenvalue decomposition of the square matrix $\mbox{F}^T\mbox{F}$
(left singular values) and form the important modes
using the dominant eigenvalues.  These modes are used as the
projection basis in the approximation given by Equation~\eqref{fapp}.
In Equation~\eqref{fapp2}, the term $ \Psi (\mbox{P}^T \Psi)^{-1} \in
\mathbb{R}^{n\times m}$ is computed once and stored.  The $d(\tau)$ is
computed using the values of the function $f(\tau)$ at $m$ points with
the indices $\wp_1,\cdots,\wp_m$ identified using the following DEIM
algorithm.\\

  \begin{tabular}{r l}
    \hline
    \hline
   DEIM& Algorithm~\cite{cs10}:\\
    \hline
    \hline

    \\
    \textbf{Input}:& The projection basis matrix $ \Psi=( \psi_{_1},\cdots,
    \psi_{_m})$ obtained by
    applying POD on a sequence \\
    & of $n_s$ function evaluations.\\
    \textbf{Output}:& The interpolation indices $\overrightarrow{\wp} = (\wp_1,\cdots,\wp_m)^T$\\
    \\
    &1: Set $[|\rho|, \wp_1] = \max\{|\psi_1|\}$\\
    &2: Set $\Psi=[\psi_1]$, $P=[e_{\wp_1}]$, and $\overrightarrow{\wp} = (\wp_1)$\\
    &3: for $k=2,...,m$ do\\
    &~~~~~- Solve $(P^T \Psi)w = P^T \psi_k$ for some $w$.\\
    &~~~~~- Compute $r = \psi_k - \Psi w$\\
    &~~~~~- Compute $[|\rho|, \wp_k] = \max\{|r|\}$\\
    &~~~~~- Set $\Psi = [\Psi ~\psi_k]$, $P = [P~ e_{\wp_k}]$,  and $\overrightarrow{\wp}=\left(
                                  \begin{array}{c}
                                    \overrightarrow{\wp} \\
                                    \wp_k \\
                                  \end{array}
                                \right)
    $\\
    & ~~end for\\
    \\
    \hline
    \hline
  \end{tabular}\\
\\
  The computational saving is due to the resulting fewer evaluations
  of $f(\tau)$. This shows the advantage of using DEIM algorithm in
  our proposed reduction method. However, applying the
    DEIM algorithm to reduce the computational cost of the nonlinear
  function requires additional computations in the offline stage,
  which will be discussed in Section~\ref{sec:GL}. Note that these
  algorithms are successful if the nonlinear functions admit low
  dimensional approximations.

%=======================

  \subsection{Generalized multiscale finite element method (GMsFEM)}
\label{sec:app}
Below we summarize the offline/online computational procedure in the
following steps:
\begin{itemize}
\item[1.]  Offline computations:
\begin{itemize}
\item 1.0. Generation coarse grid.
\item 1.1. Construction of snapshot space used to compute the offline
  space.
\item 1.2. Construction of a low dimensional space by
  performing dimension reduction in the space of local snapshots.
\end{itemize}
\item[2.] Online computations:
\begin{itemize}
\item 2.1. For each input parameter set, compute multiscale basis
  functions.
\item 2.2. Solution of a coarse-grid problem for given forcing term
  and boundary conditions.
\end{itemize}
\end{itemize}

In the offline computation, we first construct a snapshot space
$V_{\text{snap}}^{\omega_i}$.  Constructing the snapshot space may
involve solving various local problems for different choices of input
parameters or different fine-grid representations of the solution in
each coarse region.  We denote each snapshot vector (listing the
solution at each node in the domain) using a single index and create
the following matrix
$$
\Phi_{\text{snap}} = \left[ \phi_{1}^{\text{snap}}, \ldots,
  \phi_{M_{\text{snap}}}^{\text{snap}} \right],
$$
where $\phi_j^{\text{snap}}$ denotes the snapshots and
$M_{\text{snap}}$ denotes the total number of functions to keep in the
local snapshot matrix construction.

In order to construct an offline space $V_{\text{off}}$, we
reduce the dimension of the snapshot space using an
auxiliary spectral decomposition. The main objective is to use the
offline space to efficiently (and accurately) construct a set of
multiscale basis functions to be used in the online stage. More
precisely, we build a snapshot subspace that can
  approximate with sufficient accuracy any element of the original
snapshot space. The quality of the approximation is determined in the
sense defined via auxiliary bilinear forms. At the offline stage, the
bilinear forms are chosen to be \emph{parameter-independent} (through
nonlinearity), such that there is no need to reconstruct the offline
space for each $\nu$ value, where $\nu$ is assumed to be a parameter
that represents $u$ and $\mu_\kappa$ in $\kappa(x,u,\mu_\kappa)$. To
construct the offline space, we use the average of the parameters over
the coarse region $\omega_i$ in $\overline\kappa(x,\nu)$ while keeping
the spatial variations. That is, $\nu$ represents both
the average of $u$ and $\mu$.  We consider the following eigenvalue
problem in the space of snapshots,
\begin{eqnarray}
  A^{\text{off}} \Phi_k^{\text{off}} &=& \lambda_k^{\text{off}}
  S^{\text{off}}\Phi_k^{\text{off}},  \label{offeig1}
\end{eqnarray}
where
\begin{align*}
  \displaystyle A^{\text{off}}&= [a^{\text{off}}_{mn}] =
  \int_{\omega_i} \overline{\kappa}(x, \nu) \nabla
  \phi_m^{\text{snap}}
  \cdot \nabla \phi_n^{\text{snap}} = \Phi_{\text{snap}}^T \overline{A} \Phi_{\text{snap}},\\
  \displaystyle S^{\text{off}} &= [s^{\text{off}}_{mn}] =
  \int_{\omega_i} \overline{\widetilde{\kappa}}(x, \nu)
  \phi_m^{\text{snap}} \phi_n^{\text{snap}} = \Phi_{\text{snap}}^T
  \overline{S} \Phi_{\text{snap}}.
 \end{align*}
 The coefficients $\overline{\kappa}(x, \nu) $ and
 $\overline{\widetilde{\kappa}}(x, \nu)$ are parameter-averaged
 coefficients (see~\cite{egh12}). The $\overline{A}$ denotes a
 fine-scale matrix, except that parameter-averaged coefficients are
 used in its construction. The fine-scale stiffness matrix $A$ is
 constructed by integrating only on $\omega_i$
 \begin{align}
   \displaystyle A= [a_{mn}] = \int_{\omega_i} \kappa(x, u, \mu)
   \nabla \phi_m^{\text{snap}} \cdot \nabla \phi_n^{\text{snap}}.
   \label{eq:A}
 \end{align}
 To generate the offline space, we then choose the
 smallest $M_{\text{off}}$ eigenvalues from Equation~\eqref{offeig1}
 and form the corresponding eigenvectors in the respective space of
 snapshots by setting $\phi_k^{\text{off}} = \sum_j
 \Phi_{kj}^{\text{off}} \phi_j^{\text{snap}}$ (for $k=1,\ldots,
 M_{\text{off}}$), where $\Phi_{kj}^{\text{off}}$ are the coordinates
 of the vector $\Phi_{k}^{\text{off}}$. We then create the offline
 matrices
 $$
 \Phi_{\text{off}} = \left[ \phi_{1}^{\text{off}}, \ldots,
   \phi_{M_{\text{off}}}^{\text{off}} \right]
$$
to be used in the online space construction.

The online coarse space is used within the finite element framework to
solve the original global problem, where continuous Galerkin
multiscale basis functions are used to compute the global solution. In
particular, we seek a subspace of the respective offline space such
that it can approximate well any element of the offline space in an
appropriate metric. At the online stage, the bilinear forms are chosen
to be \emph{parameter dependent}.  The following eigenvalue problems
are posed in the reduced offline space:
\begin{eqnarray}
  A^{\text{on}}(\nu) \Phi_k^{\text{on}} &=& \lambda_k^{\text{on}}
  S^{\text{on}}(\nu)\Phi_k^{\text{on}},  \label{oneig1}
\end{eqnarray}
where
\begin{equation*}
  \displaystyle A^{\text{on}}(\nu) = [a^{\text{on}}(\nu)_{mn}] =
  \int_{\omega_i} \kappa(x, \nu) \nabla \phi_m^{\text{off}} \cdot
  \nabla \phi_n^{\text{off}} = \Phi_{\text{off}}^T A(\nu)
  \Phi_{\text{off}},
 \end{equation*}
 
 \begin{equation*}
  {\displaystyle S^{\text{on}}(\nu) = [s^{\text{on}}(\nu)_{mn}] =
  \int_{\omega_i} \widetilde{\kappa}(x, \nu)  \phi_m^{\text{off}} 
  \phi_n^{\text{off}} = \Phi_{\text{off}}^T S(\nu)
  \Phi_{\text{off}},}
 \end{equation*}
 
 and $\kappa(x, \nu)$ and 
   $\widetilde{\kappa}(x, \nu)$ are now parameter dependent.  To generate the
 online space, we then choose the smallest $M_{\text{on}}$ eigenvalues
 from \eqref{oneig1} and form the corresponding eigenvectors in the
 offline space by setting $\phi_k^{\text{on}} = \sum_j
 \Phi_{kj}^{\text{on}} \phi_j^{\text{off}}$ (for $k=1,\ldots,
 M_{\text{on}}$), where $\Phi_{kj}^{\text{on}}$ are the coordinates of
 the vector $\Phi_{k}^{\text{on}}$.  If $\kappa(x,u)=k_0(x) b(u)$,
 then one can use the parameter-independent case of GMsFEM. In this
 case, there is no need to construct an online space (i.e., the
 online space is the same as the offline space). From now on, we
 denote the online space basis functions by $\phi_i$.

\section{Global-Local Nonlinear Model Reduction}
\label {sec:globallocal}

%=======================
\subsection{Local multiscale model reduction }
\label{sec:Newton}
The finite element discretization of~\eqref{parabolic} yields a system
of ordinary differential equations given by
\begin{eqnarray}\label{sode}
  \mbox{M} \dot{\mbox{U}} + \mbox{F}(\mbox{U}) = \mbox{H},
\end{eqnarray}
where
$$
\mbox{U} = \left(
               \begin{array}{cccc}
                 u_1 & u_2 & \cdots & u_{N_f} \\
               \end{array}
             \right)
$$
is the vector collecting the pressure values at all nodes in the local
domain and $ \mbox{H}$ is the right-hand-side vector obtained by
discretization.  Using the offline basis functions, we can write (in
discrete form)
\begin{equation}\label{eqjuan:kappa}
  \kappa(x,u,\mu)= \sum_{q=1}^Q \kappa_q(x) b_q(u,\mu).
\end{equation}
This results in
\begin{equation*}
  \mbox{F}(\mbox{U},\mu) = \sum_{q=1}^Q\mbox{A}_q
  \Lambda_1^q(\mbox{U},\mu) \mbox{U},
\end{equation*}
 where we have
\begin{align*}
  \mbox{A}_q:= [\mbox{a}_{ij}^q]&=\int_\Omega \kappa_q \nabla \phi^0_i
  \cdot \nabla \phi^0_j , \;\;\; \mbox{M}:=[\mbox{m}_{ij}]=\int_\Omega
  \phi^0_i \phi^0_j\;,
  \;\;\; \mbox{H}:=[\mbox{h}_{i}]=\int_\Omega  \phi^0_i h \;, \\
  \Lambda_1^q(\mbox{U},\mu) &= \mbox{diag}\left(
               \begin{array}{cccc}
                 b_q(u_1,\mu) & b_q(u_2,\mu) & \cdots &
                 b_q(u_{N_f},\mu) \\
               \end{array}
             \right),
\end{align*}
and $\phi^0_i$ are piecewise linear basis functions defined on a fine
triangulation of $\Omega$.

Employing the backward Euler scheme for the time marching process, we
obtain
\begin{eqnarray}
  \mbox{U}^{n+1}+\Delta t \;\mbox{M}^{-1}\mbox{F}(\mbox{U}^{n+1})=
  \mbox{U}^{n} + \Delta t \;\mbox{M}^{-1} \mbox{H},
\end{eqnarray}
where $\Delta t$ is the time-step size and the superscript $n$ refers
to the temporal level of the solution. The residue is
  defined as:
\begin{eqnarray}
  \mbox{R}(\mbox{U}^{n+1}) = \mbox{U}^{n+1}- \mbox{U}^{n} + \Delta t
  \;\mbox{M}^{-1}\mbox{F}(\mbox{U}^{n+1}) - \Delta t \;\mbox{M}^{-1}
  \mbox{H}
\end{eqnarray}
with derivative (Jacobian)
\begin{align}
\mbox{J}(\mbox{U}^{n+1}) = D\mbox{R}(\mbox{U}^{n+1}) &= I + \Delta t
\;\mbox{M}^{-1} D \mbox{F}(\mbox{U}^{n+1}) \nonumber \\
&= I + \sum_{q=1}^Q \Delta t \;\mbox{M}^{-1}
\mbox{A}_q\Lambda_1^q(\mbox{U}^{n+1}) +\sum_{q=1}^Q \Delta t
\;\mbox{M}^{-1} \mbox{A}_q\Lambda_2^q(\mbox{U}^{n+1}),
\end{align}
where
\begin{equation*}
\Lambda_2^q(\mbox{U},\mu) = \mbox{diag}\left(
               \begin{array}{cccc}
                u_1  \frac{\partial b_q(u_1,\mu)}{\partial u} &
                u_2  \frac{\partial b_q(u_2,\mu)}{\partial u} &
                   \cdots &
                 u_{N_f}   \frac{\partial b_q(u_{N_f},\mu)}{\partial u}
               \end{array}
             \right),
\end{equation*}
and $D$ is the multi-variate gradient operator defined as
$\left[D\mbox{R}(\mbox{U})\right]_{ij}=\partial {R}_i/\partial {U}_j$.
%$$
%\Lambda_2(\mbox{U}) = \mbox{diag}\left(
%               \begin{array}{cccc}
%                 \mu u_1 e^{\mu u_1} & \mu u_2 e^{\mu u_2} & \cdots & \mu u_{N_f} e^{\mu u_{N_f}} \\
%               \end{array}
%             \right).
%$$
The scheme involves, at each time step, the following iterations
\begin{align}
  \mbox{J}(\mbox{U}^{n+1}_{(k)}) \Delta \mbox{U}^{n+1}_{(k)}&=
  -\Big(\mbox{U}^{n+1}_{(k)}- \mbox{U}^{n} + \Delta t
  \;\mbox{M}^{-1}\mbox{F}(\mbox{U}^{n+1}_{(k)}) - \Delta t
  \;\mbox{M}^{-1} \mbox{H}\Big) \nonumber \\
  \mbox{U}^{n+1}_{(k+1)}&=\mbox{U}^{n+1}_{(k)} + \Delta
  \mbox{U}^{n+1}_{(k)}, \nonumber
\end{align}
where the initial guess is $\mbox{U}^{n+1}_{(0)}=\mbox{U}^{n}\nonumber
$ and $k$ is the iteration counter. The above iterations are
repeatedly applied until $\parallel \Delta
\mbox{U}^{n+1}_{(k)} \parallel$ is less than a specific tolerance.

In our simulations, we use $Q=1$ in~\eqref{eqjuan:kappa} as our focus
is on localized multiscale interpolation of nonlinear functionals that
arise in discretization of multiscale PDEs.  With this choice, we do
not need to compute the online multiscale space (i.e., the online
space is the same as the offline space).

We use the solution expansion (i.e., $u=\Phi z$)
and employ the multiscale framework to obtain a set of $N_c$ ordinary
differential equations that constitute a reduced-order model; that is,
\begin{equation}\label{rom11}
  \dot{z} = - (\Phi^T\mbox{M}\Phi)^{-1} \Phi^T \mbox{F} (\Phi z) +
  (\Phi^T\mbox{M}\Phi)^{-1} \Phi^T \mbox{H}.
\end{equation}
Thus, the original problem with $N_f$ degrees of freedom is reduced to
a dynamical system with $N_c$ dimensions where $N_c \ll N_f$.

The nonlinear term $(\Phi^T\mbox{M}\Phi)^{-1} \Phi^T \mbox{F} (\Phi
z)$ in the reduced-order model, given by Equation (\ref{rom11}), has a
computational complexity that depends on the dimension of the full
system $N_f$.  As such, solving the reduced
system still requires extensive computational resources and time. To
reduce this computational requirement, we use multiscale DEIM as
described in the previous section.
%In this case, computational savings
%can be obtained in a forward run of the nonlinear model.

To solve the reduced system, we employ the backward Euler
scheme; that is,
\begin{eqnarray}
z^{n+1}+\Delta t \;\widetilde{\mbox{M}}^{-1}
\widetilde{\mbox{F}}(z^{n+1})= z^{n} + \Delta t
\;\widetilde{\mbox{M}}^{-1} \widetilde{\mbox{H}},
\end{eqnarray}
where $\widetilde{\mbox{M}}= \Phi^T\mbox{M}\Phi$,
$\widetilde{\mbox{F}}(z)=\Phi^T\mbox{F}(\Phi z)$, and
$\widetilde{\mbox{H}}=\Phi^T \mbox{H}$. We let
\begin{eqnarray}
\label{eq:msres}
\widetilde{\mbox{R}}(z^{n+1}) = z^{n+1}- z^{n} +\Delta t
\;\widetilde{\mbox{M}}^{-1} \widetilde{\mbox{F}}(z^{n+1}) - \Delta t
\;\widetilde{\mbox{M}}^{-1} \widetilde{\mbox{H}}
\end{eqnarray}
with derivative
\begin{align}
\widetilde{\mbox{J}}(z^{n+1})=D\widetilde{\mbox{R}}(z^{n+1}) &= I +
\Delta t \;\widetilde{\mbox{M}}^{-1} D \widetilde{\mbox{F}}(z^{n+1})
\nonumber \\
&= I + \sum_{q=1}^Q\Delta t \;\widetilde{\mbox{M}}^{-1}
\Phi^T\mbox{A}_q\Lambda_1^q(\Phi z^{n+1})\Phi +\sum_{q=1}^Q \Delta t
\;\widetilde{\mbox{M}}^{-1} \Phi^T\mbox{A}_q\Lambda_2^q(\Phi
z^{n+1})\Phi. \nonumber
\end{align}
The scheme involves, at each time step, the following iterations
\begin{align}\label{scheme}
  \widetilde{\mbox{J}}(z^{n+1}_{(k)}) \Delta z^{n+1}_{(k)}&=
  -\Big(z^{n+1}_{(k)}- z^{n} + \Delta t
  \;\widetilde{\mbox{M}}^{-1}\widetilde{\mbox{F}}(z^{n+1}_{(k)}) -
  \Delta t \;\widetilde{\mbox{M}}^{-1} \widetilde{\mbox{H}} \Big) \\
  z^{n+1}_{(k+1)}&=z^{n+1}_{(k)} + \Delta z^{n+1}_{(k)},
\end{align}
where the initial guess is $z^{n+1}_{(0)}=z^{n}\nonumber $. The above
iterations are repeated until $\parallel \Delta
z^{n+1}_{(k)} \parallel$ is less than a specific
tolerance. Furthermore, we use multiscale DEIM to
approximate the nonlinear functions that appear in the residual
$\widetilde{\mbox{R}}$ and the Jacobian $\widetilde{\mbox{J}}$ to
reduce the number of function evaluations.

\subsection{Global-local nonlinear model reduction approach}
\label{sec:GL}
We denote the offline parameters by $\theta^{\text{off}}$ which
include samples of the right-hand side $h(x)$ denoted by
$h_i^{\text{off}}$, samples of $\mu$ denoted by $\mu_i^{\text{off}}$,
and samples of initial conditions denoted by
$U_{0,i}^{\text{off}}$. Similarly, the online parameter set is denoted
by $\theta^{\text{on}}$ and includes the online source term
$h^{\text{on}}$, the online $\mu$ ($\mu^{\text{on}}$), and the online
initial conditions $U_0^{\text{on}}$. We follow a global-local
nonlinear model reduction approach that includes the following steps:
\begin{itemize}
\item \textbf{Offline Stage}
\\
The offline stage includes the following steps:
\begin{itemize}

\item  Consider the offline parameters set $\theta^{\text{off}}=\{\theta_i^{\text{off}}\}=\{h^{\text{off}}_i,\mu^{\text{off}}_i, U_{0,i}^{\text{off}}\}$.

\item  Use $\theta^{\text{off}}_i$ to define the fine-scale stiffness and mass matrices, source terms and multiscale basis functions.
\item  Compute the local snapshots of the nonlinear functions and use DEIM algorithm, as described in the previous section, to set the local DEIM basis functions and local DEIM points ($L_0^{local}$).
\item Generate snapshots of the coarse-grid solutions using local DEIM.

\item Record $N_t$ instantaneous solutions (usually referred as snapshots) using coarse-grid approximations from the above step and collect them in a snapshot matrix as:
\begin{eqnarray}
\textbf{Z}^{N_t}=\{\textbf{Z}_1,\textbf{Z}_2,\textbf{Z}_3,\cdots, \textbf{Z}_{N_t}\}
\end{eqnarray}
where $N_t$ is the number of snapshots and $N_c$ is the size of the column vectors $\textbf{Z}_{i}$.

\item Compute the POD modes and use these modes to approximate the solution field on the coarse grid. As such, we assume an expansion in terms of the modes $\psi_i$; that is, we let
\begin{equation}
z(x,t) \approx \tilde{z}(x,t) = \sum_{i=1}^{N_r} \alpha_i(t) \psi_i(x)\label{utilda}
\end{equation}
or in a matrix form
\begin{equation}
\mathbf{Z}^n \approx \tilde{\mathbf{Z}}^n = \Psi \alpha^n \label{utildam}
\end{equation}
where $\Psi=\left(
              \begin{array}{ccc}
                \psi_1 & \cdots & \psi_{N_r} \\
              \end{array}
            \right)$.

\end{itemize}

\item \textbf{ Online Stage}
\\
The online stage includes the following steps:
\begin{itemize}
\item Given online $\theta ^{\text{on}}=\{h^{\text{on}}, \mu^{\text{on}},U_0^{\text{on}}\}$
\item Use the solution expansion given by~\eqref{utilda}
  and project the governing equation of the coarse-scale problem onto
  the space formed by the modes to obtain a set of $N_r$ ordinary
  differential equations that constitute a reduced-order model; that
  is,
\begin{equation}\label{rom}
  \dot{\alpha} = - (\underbrace{\Psi^T\Phi^T\mathbf{M}\Phi\Psi}_{N_r\times N_r})^{-1}\underbrace{\Psi^T }_{N_r\times N_c}\underbrace{\Phi^T}_{N_c\times N_f} \underbrace{\mathbf{F}(\Phi\Psi\alpha)}_{N_f \times1}+ (\Psi^T\Phi^T\mathbf{M}\Phi\Psi)^{-1}\Psi^T \Phi^T \mathbf{H}.
\end{equation}
\item Employ Newton's method to solve the above reduced system. The
  Newton scheme involves at each time step the following iteration. We
  need to solve the linear system
\begin{align}\label{scheme1}
  \widehat{\mbox{J}}(\alpha^{n+1}_{(k)}) \Delta \alpha^{n+1}_{(k)}&=
  -\Big(\alpha^{n+1}_{(k)}- \alpha^{n} + \Delta t
  \;\widehat{\mbox{M}}^{-1}\widehat{\mbox{F}}(\alpha^{n+1}_{(k)}) -
  \Delta t \;\widehat{\mbox{M}}^{-1} \widehat{\mbox{H}} \Big)
\end{align}
where
$$
\widehat{\mbox{M}}=\Psi^T\Phi^T\mathbf{M}\Phi\Psi=\Psi^T\tilde{\mathbf{M}}\Psi,~~~~~
\widehat{\mbox{H}}=\Psi^T\Phi^T\mathbf{H}=\Psi^T\tilde{\mathbf{H}},~~~~~
\widehat{\mbox{F}}=\Psi^T\Phi^T\mathbf{F}=\Psi^T\tilde{\mathbf{F}}.$$
Then
\[
\alpha^{n+1}_{(k+1)}=\alpha^{n+1}_{(k)} -(
\widehat{\mbox{J}}(\alpha^{n+1}_{(k)}) )^{-1} \Big(\alpha^{n+1}_{(k)}-
\alpha^{n} + \Delta t
\;\widehat{\mbox{M}}^{-1}\widehat{\mbox{F}}(\alpha^{n+1}_{(k)}) -
\Delta t \;\widehat{\mbox{M}}^{-1} \widehat{\mbox{H}} \Big).
\]
Thus, the original problem with $N_f$ degrees of freedom
is reduced to a dynamical system with $N_r$ dimensions where $N_r \ll
N_c\ll N_f$.

\item Use global DEIM to approximate the nonlinear functions that
  appear in the residual and Jacobian. To do so, we write the
  nonlinear function ${\mathbf{F}(\Phi\Psi\alpha)}$ in
  Equation~\eqref{rom} as
 \begin{equation}\label{fapp3}
  \mathbf{F}(\Phi\Psi\alpha) \approx  \Psi^*  {d},
\end{equation}
where $\Psi^*=[\psi_1^*,...,\psi_{L_0^{\text{global}}}^*]$ is the
matrix of the global DEIM basis functions
$\{\psi_i^*\}_{i=1}^{L_0^{\text{global}}}$. These functions are
constructed using the snapshots of the nonlinear function
$\mathbf{F}(\Phi z)$ computed offline and employ the POD technique to
select the most energetic modes (see Section~\ref{SEC:DEIM}). The
coefficient vector ${d}$ is computed using the values of the function
$\mathbf{F}$ at $L_0^{\text{global}}$ global points.
\item Use the solution expansion given by~\eqref{utilda} in terms of
  POD modes to approximate the coarse-scale solution and then use the
  operator matrix $\Phi$ to downscale the approximate solution and
  evaluate the flow field on the fine grid.

\end{itemize}
\end{itemize}
\section{Numerical Results}
\label{sec:num}
In this section, we use representative numerical
examples to illustrate the applicability of the proposed global-local
nonlinear model reduction approach for solving nonlinear multiscale
partial differential equations. Before presenting the individual
examples, we describe the computational domain used in constructing
the GMsFEM basis functions. This computation is performed during the
offline stage.  We discretize with linear finite elements a nonlinear
PDE posed on the computational domain $D=[0,1]\times[0,1]$.  For
constructing the coarse grid, we divide $[0,1]\times[0,1]$ into
$10\times 10$ squares. Each square is divided further into
$10\times10$ squares each of which is divided into two triangles.
Thus, the mesh size is $1/100$ for the fine mesh and
  $1/10$ for the coarse one. The fine-scale finite element vectors
introduced in this section are defined on this fine grid. The
fine-grid representation of a coarse-scale vector $z$ is given by
$\Phi z$, which is a fine-grid vector.

\begin{figure}[htb]
  \centering
  \includegraphics [width=0.47\textwidth]{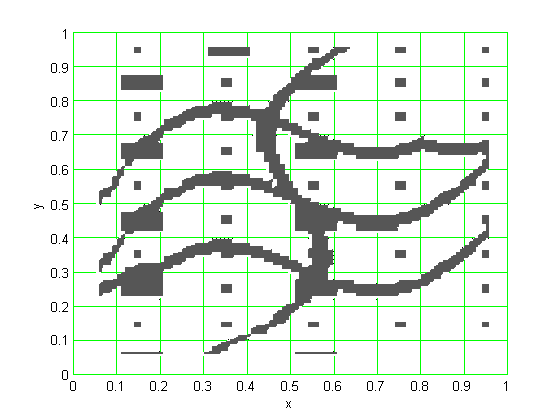}
  \caption{\label{perm1}{Permeability field that model high
      conductivity channels within a homogeneous domain. The minimum
      (background) conductivity is taken to be $\kappa_{min}=1$, and
      the high conductivity (gray regions) with value of
      $\kappa_{max}=\eta~(\eta=10^6)$}}
\end{figure}

In the following numerical examples, we consider~\eqref{parabolic}
with specified boundary and initial conditions, where the permeability
coefficient and the forcing term are given by
\begin{center}
  $\kappa(x;u,\mu)=\kappa_q(x)b_q(u,\mu)$ \ and \ $h(x)=1+\sin(2\pi
  x_1)\sin(2\pi x_2)$.
\end{center}
Here, $\kappa_q$ represents the permeability field with
high-conductivity channels as shown in Figure~\ref{perm1} and
$b_q(u,\mu)$ is defined later for each example. We use the GMsFEM
along with the Newton method to
discretize~\eqref{parabolic}. Furthermore, we employ the local
multiscale DEIM in the offline stage and the global multiscale DEIM in
the online stage to approximate the nonlinear functions that arise in
the residual and the Jacobian.

Using the fine-scale stiffness matrix $\mbox{A}$ that corresponds
to~\eqref{parabolic}, as defined in~\eqref{eq:A}, we
  introduce the relative energy error as
\begin{eqnarray}\label{Err2}
  \|E\|_{\mbox{A}} = \sqrt{\frac{(\mbox{U} - \widetilde{\mbox{U}})^T
      {\mbox{A}} (\mbox{U} - \widetilde{\mbox{U}})}{\mbox{U}^T
      {\mbox{A}} \mbox{U}}}.
\end{eqnarray}
Moreover, we define $w_0$ to be the solution of the problem
\begin{equation} \label{elliptic}
 - \nabla \cdot \left( \kappa_q(x) \nabla w_0 \right) = h(x)  \quad \text{in} \quad D,
\end{equation}
to use it in the following examples as our initial guess.  In the
following, we show:
%In general, these examples we show the following:
\begin{itemize}
\item In the first example, we compare the approximate solution of the
  reduced system obtained by applying the global-local
  approach against the solution of the original system with full
  dimension ($N_f$) and show the reduction we achieve in terms of the
  computational cost.
\item In the second example, we show the variations of the error as we
  increase the number of local DEIM points, $L_0^{local}$, and global
  DEIM points, $L_0^{global}$ for one selection of the parameter
  $\mu$.
\item In the third example, we show the effect of using several
  offline parameters to improve the reduced-order solutions. As such,
  we use two offline values of the parameter $\mu$ and solve an online
  problem for a different value of $\mu$.
\item In the fourth example, we use two offline values of $\mu$ and
  show the variations of the errors as we increase the number of local
  and global points.
\item {Random values of the parameter $\mu$ with a probability distribution are used in the fifth example. We demonstrate the applicability of our approach in this setup.}
\end{itemize}

\subsection{Single Offline Parameter}
\begin{exmp}
\label{ex1}
We consider~\eqref{parabolic} along with the following
offline and online parameters
\[
\theta^{\text{off}}=
\begin{cases}
h^{\text{off}}= &  1+\sin(2\pi x_1)\sin(2\pi x_2), \\
\mu^{\text{off}}= & 10,\\
 U_0^{\text{off}}= & w_0,
         \end{cases}
          \]\\
  \[
\theta^{\text{on}}=
\begin{cases}
  h^{\text{on}}= &  1+\sin(2\pi x_1)\sin(2\pi x_2), \\
  \mu^{\text{on}}= & 40,\\
  U_0^{\text{on}}= & w_0*0.5.
\end{cases}
\]
where the nonlinear function $b_q$ is defined as $ b_q(u,\mu)=e^{\mu u
}$. Here, the source term does not need to be fixed for the method to
work as we see below. We employ GMsFEM for the spatial discretization
and the backward Euler method for time advancing as described in
Section \ref{sec:Newton}. Furthermore, we follow the steps given in
Section \ref{sec:GL} using three DEIM points ($L_0^{local}=3$) per
coarse region to approximate $b_q$ in the offline stage. After
generating the snapshots of the coarse-grid solutions using local
DEIM, we compute the multiscale POD modes that are used in the online
problem. We use $L_0^{global}=5$ in the online stage to approximate
$b_q$ globally and then use the generated POD modes to approximate the
coarse-scale solution. In Figure \ref{sol}, we compare the approximate
solution obtained from the global-local nonlinear model
reduction approach with the solution of the original system without
using the DEIM technique to approximate the nonlinear function. A good
approximation is observed  in  this figure, which demonstrates the
capability of global-local nonlinear model reduction to reproduce
accurately the fully resolved solution of a nonlinear PDE.

We have also considered a permeability field that is obtained by
rotating the permeability field $\kappa_q$ in Figure~\ref{perm1} such
that the three long channels are in the vertical direction.  Our
numerical results show similar accuracy and computational cost
compared to the previous case (see Figure~\ref{perm1}). In general, we
expect non-homogeneous boundary conditions to affect the
numerical results.

The approximate solution shown in Figure~\ref{approxsol} is obtained
using only two POD modes. As expected, increasing the number of POD
modes used in the online stage yields a better
approximation. That is, the error decreases as we increase the number
of POD modes used as shown in Figure~\ref{POD}. The error using two
POD modes decreases slightly from $12\%$ (at steady state) to $11.5\%$
when using three POD modes. The decreasing trend is steeper when
considering more POD modes. For instance, the use of $5$ modes yields
an error of $4.5\%$.

\begin{figure}[htp]
  \begin{center}
    \subfigure[Reference
    Solution]{\includegraphics[width=0.47\textwidth]{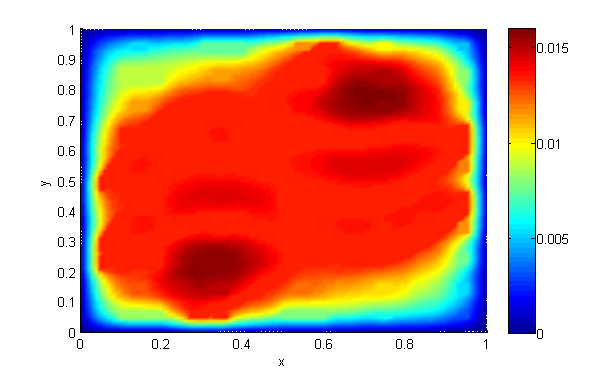}\label{refsol}}
    \subfigure[Approximate Solution
    ]{\includegraphics[width=0.47\textwidth]{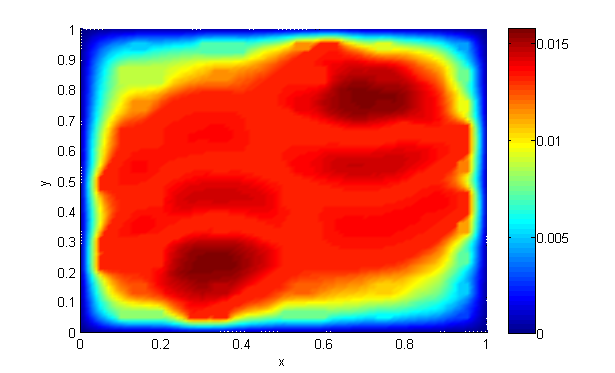}\label{approxsol}}
  \end{center}

  \caption{\label{sol}{Comparison between reference
      solution of the fine-scale problem with that obtained from the
      global-local multiscale approach.}}

\end{figure}

\begin{figure}[htp]
  \begin{center}
    \includegraphics[width=0.47\textwidth]{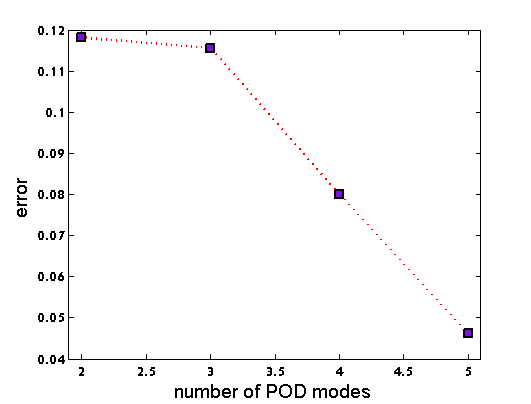}
  \end{center}
  \caption{\label{POD}{Variations of the solution error
      with the number of POD modes.}}
\end{figure}

In order to illustrate the computational savings, we compute the time
for solving the system of ordinary differential equations given
in~\eqref{sode} with and without using the proposed method. We denote
the time for solving the full system by $T_{fine}$ and the time for
solving the reduced system using global-local nonlinear model
reduction by $T_{GL}$. Then, the percentage of the simulation time is
given by
 \begin{equation}\label{rtsp}
  R=\frac{T_{GL} }{T_{fine}}*100.
\end{equation}
We compute $R$ with respect to different number of DEIM points and POD
modes and present the results in Tables~\ref{tabl1a} and~\ref{tabl1b},
respectively. In Table~\ref{tabl1a}, the first column shows the number
of local DEIM points ($L_0^{local}$), the second column represents the
number of global DEIM points ($L_0^{global}$), and the third column
illustrates the percentage of the simulation time. Here two POD modes
are used. As $L_0^{local}$ and/or $L_0^{global}$ increase, the
percentage decreases accordingly. For example, $R$ decreases from
$3.7832~\%$ to $3.3741~\%$ by increasing $L_0^{global}$ from two to
three, and to $3.2093~\%$ by increasing both $L_0^{local}$ and
$L_0^{global}$ from two to three. Decreasing $R$ means that $T_{GL}$,
time for solving the reduced system, decreases as we
  increase the number of DEIM points. Therefore, increasing the
number of local and global DEIM points may speed up the
simulation in addition to improving the accuracy as we see in the next
example. In Table~\ref{tabl1b}, the numbers of POD modes used for the
global reduction are listed in the first column and the corresponding
values of $R$ are shown in the second column. In this case, we keep the number
of local and global DEIM points constant and equal to two and three,
respectively. Now, increasing the number of POD modes inversely
affects the simulation speed-up. That is, increasing the number of POD
modes increases the value of $R$ which means $T_{GL}$ is increasing
and hence the speed-up of our simulation is decreasing. For example,
$R$ increases from $3.3741~\%$ when we use two POD modes to
$4.0387~\%$ with three POD modes and keeps increasing as
we increase the number of POD modes to be $6.1414~\%$ with five POD
modes.  Although, increasing the number of POD modes slows down the
simulation, it improves the accuracy of the approximate solution (see
Figure~\ref{POD}). However, the following examples show the capability
in terms of the accuracy of this method when using two POD modes for
the global reduction.
\begin{table}[htp]
  \centering
  \begin{tabular}{|c|c|c|}
    \hline

    $L_0^{local}$ & $L_0^{global}$   &  $R(\%)$ \\
    \hline
    $2$ & $ 2$& $3.7832$\\
    \hline
    $2$& $3$ & $ 3.3741$\\
    \hline
    $3$& $3$ & $3.2093$\\
    \hline
  \end{tabular}
  \caption{ Variation of the percentage of the simulation time corresponding to different number  of local and global DEIM points. Here we use two POD modes.}
  \label{tabl1a}
\end{table}

\begin{table}[htp]
\centering
\begin{tabular}{|c|c|}
    \hline

     POD modes  & $R(\%)$  \\
    \hline
    $2$ & $3.3741$\\
    \hline
    $3$    &   $ 4.0387$\\
    \hline
    $4$    &   $4.9158$\\
    \hline
    $5$    &  $6.1414$\\
    \hline
  \end{tabular}
  \caption{ Variation of the percentage of the simulation time corresponding to different number of POD
modes. Here we use $L_0^{local}=2$ and $L_0^{global}=3$}
\label{tabl1b}
\end{table}

\end{exmp}
\begin{exmp}
\label{ex2}
\*
In this example, we use different numbers of local and global DEIM points, $L_0^{local}=\{1, 2, 3\}$ and
$L_0^{global}=\{1, 2, 3\}$, to investigate how these numbers affect the error. As in Example \ref{ex1}, we
consider $ b_q(u,\mu)=e^{\mu u }$ and the following offline and online parameters:
\[
\theta^{\text{off}}=
\begin{cases}
h^{\text{off}}= &  1+\sin(2\pi x_1)\sin(2\pi x_2), \\
\mu^{\text{off}}= & 10,\\
 U_0^{\text{off}}= & w_0,
         \end{cases}
          \]\\
  \[
\theta^{\text{on}}=
\begin{cases}
h^{\text{on}}= &  1+\sin(2\pi x_1)\sin(2\pi x_2), \\
\mu^{\text{on}}= & 40,\\
 U_0^{\text{on}}= & 0.5 w_0.
         \end{cases}
          \]
\begin{figure}[t]
  \begin{center}
    \subfigure[Variations of global DEIM points]{\includegraphics[width=0.47\textwidth]{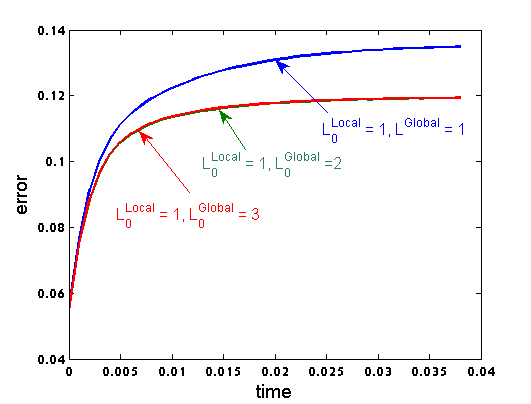}\label{2a}}
    \subfigure[Variations of local DEIM points]{\includegraphics[width=0.47\textwidth]{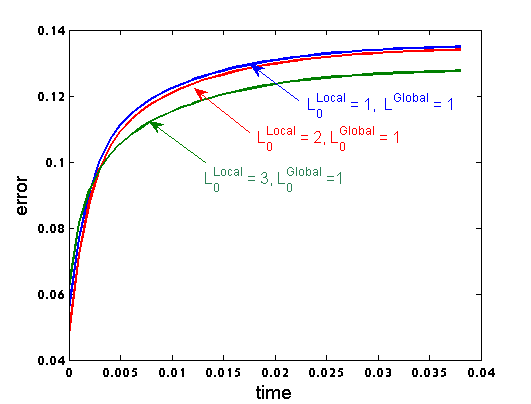}\label{2b}}
    \subfigure[Variations of both global and local DEIM points]{\includegraphics[width=0.47\textwidth]{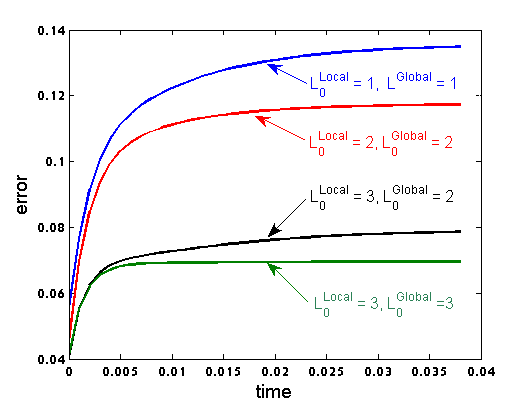}\label{2c}}

  \end{center}
  \caption{\label{errorDEIMpoints}{Effect of the number of local and global DEIM points on the approximate solution accuracy.}}

\end{figure}

In Figure \ref{2a}, we plot the transient variations of the error
while using different numbers of global DEIM points for a
  fixed number of local DEIM points equal to one. Increasing the
number of global DEIM points from one to three results in
a decrease in the error from $13\%$ to $11\%$ (at steady
state). Further increases in the number of global DEIM points does not
yield any improvement in the total error. This is due to the dominance
of the local error. Figure \ref{2b} shows the decreasing trend of the
error as we increase the number of local DEIM points. In Figure
\ref{2c}, we show the variations of the error with increasing the
number of both local and global DEIM points. Increasing the number of
DEIM points enables smaller error and then improves the solution
accuracy. These examples show that the number of local and global DEIM
points need to be chosen carefully to balance the local and global
errors.
\\

\end{exmp}

\subsection{Multiple Offline Parameters}
\begin{exmp}
\label{ex3}
In this example, we define the nonlinear function as
$b_q(u,\mu)=e^{\mu(0.9+u)}$ and use $\mu^{\text{off}}_1=2$ and
$\mu^{\text{off}}_2=5$, separately, in the offline problem to compute
POD modes and DEIM points. We then combine these to use the total
number of POD modes in the online problem with a different online
value of $\mu$ ($\mu^{\text{on}}=3$). In this example we keep the
number of local and global DEIM points constant and equal to three
(i.e., $L^{local}_0=L^{global}_0=3$).  Furthermore, we use different
online initial conditions and source term. The following system
parameters are considered.

\[
\theta^{\text{off}}=
\begin{cases}
h^{\text{off}}= &  1+\sin(2\pi x_1)\sin(2\pi x_2), \\
\mu^{\text{off}}_1= & 2,\\
\mu^{\text{off}}_2= & 5,\\
 U_0^{\text{off}}= & w_0.
         \end{cases}
          \]\\
  \[
\theta^{\text{on}}=
\begin{cases}
h^{\text{on}}= &  1+\sin(4\pi x_1)\sin(4\pi x_2), \\
\mu^{\text{on}}= & 3,\\
 U_0^{\text{on}}= & 0.
         \end{cases}
          \]

\begin{figure}[t]
  \begin{center}
    \includegraphics[width=0.47\textwidth]{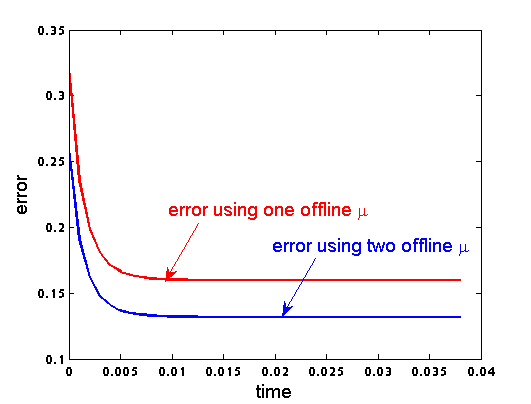}\label{21}
       \end{center}
  \caption{\label{diffmu1}{Transient variations of the error (using different offline values of the parameter $\mu$).}}
\end{figure}

We show in Figure \ref{diffmu1} that the error decreases when combining two cases that correspond to different values of offline $\mu$. For instance, the error when considering only one offline case is about $16\%$ and it goes down to $13\%$ when combining two cases with two different values of offline $\mu$. Hence, using multiple parameter values in the offline stage improves the method's accuracy independently of the online parameters.
\end{exmp}
\begin{exmp}
\label{ex4}
 Next, we consider the following parameters
 \[
\theta^{\text{off}}=
\begin{cases}
h^{\text{off}}= &  1+\sin(2\pi x_1)\sin(2\pi x_2), \\
\mu^{\text{off}}_1= & 10,\\
\mu^{\text{off}}_2= & 40,\\
 U_0^{\text{off}}= & w_0,
         \end{cases}
          \]\\
  \[
\theta^{\text{on}}=
\begin{cases}
h^{\text{on}}= &  1+\sin(2\pi x_1)\sin(2\pi x_2), \\
\mu^{\text{on}}= & 24,\\
 U_0^{\text{on}}= & 0.
         \end{cases}
          \]

\noindent
 and the nonlinear function $b_q(u,\mu)=e^{\mu u}$. In this case, we use two offline values of $\mu$ while considering different numbers of local and global DEIM points.
\begin{figure}[t]
  \begin{center}
    \subfigure[Variations of global DEIM points]{\includegraphics[width=0.47\textwidth]{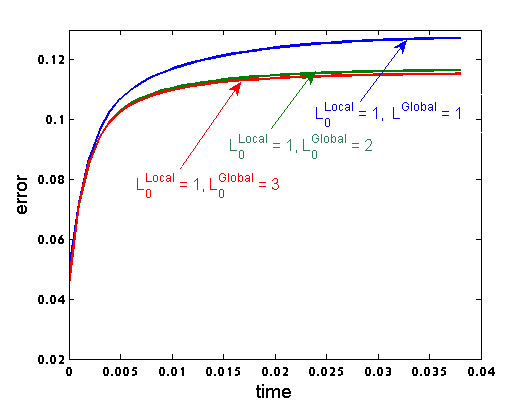}\label{3a}}
    \subfigure[Variations of local DEIM points]{\includegraphics[width=0.47\textwidth]{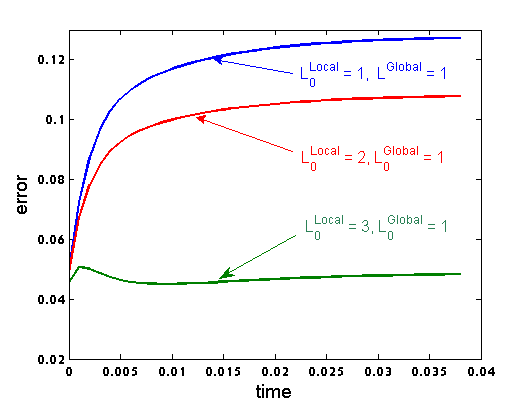}\label{3b}}
    \subfigure[Variations of both global and local DEIM points]{\includegraphics[width=0.47\textwidth]{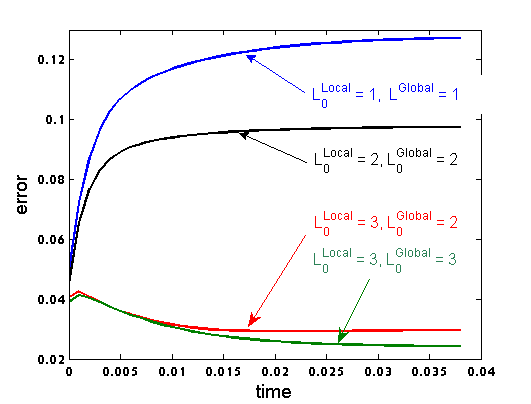}\label{3c}}
  \end{center}
  \caption{\label{ex3a}{Effect of the number of local and global DEIM points on the approximate
solution accuracy (using two offline $\mu$).}}

\end{figure}
The effect of the number of local and global DEIM points on the error
between the reference and approximate solutions when combining two
cases that correspond to two different values of $\mu$ is shown in
Figure \ref{ex3a}. Similar trends to those of Example \ref{ex2} are
observed. Increasing both local and global DEIM points improves the
approximation to the solution. For instance, the error reduces from
about $13\%$ when using a local and a global DEIM point to $2\%$ when
using three local and global DEIM points. The error reduction in this
case (when we use two offline $\mu$) is bigger than the one we
obtained when only using one offline $\mu$ value where the error
decreased from $13\%$ to $7\%$ (see Figure \ref{2c}). We conclude that using two offline
$\mu$ values and increasing number of local and global DEIM points
yields a better approximation. Therefore, choosing the number of local
and global DEIM points and the offline parameter values are the main
factors to achieve high accuracy in the  proposed method.
\end{exmp}

\begin{exmp}
\label{ex5}

In this example we consider the case with random values of the
parameter $\mu$ that has a normal distribution with the mean $25$ and
variance $4$.  As in Example \ref{ex3}, we use different values of the
offline parameter $\mu^{\text{off}}=\{10,25,39\}$, and compute the POD
and DEIM modes. Further, we combine these modes to get the global POD
and DEIM modes that we use in the online problem. In the
online problem, we take uncorrelated random values of
$\mu^{\text{on}}$ drawn from the above probability distribution.  We
rapidly compute the approximate solution and 
evaluate the relative error corresponding
to each value of $\mu^{\text{on}}$. Comparing the mean solutions of
the fully-resolved model and the reduced model demonstrates the
capability of the proposed method when random values of the parameter
is employed in the nonlinear functional. 
Furthermore, we observe a good accuracy
as  shown from the error plotting in Figure \ref{meanerror}.

%
%\begin{figure}[htp]
%  \begin{center}
%    \subfigure[Mean Reference Solution]{\includegraphics[width=0.47\textwidth]{refsol}\label{MeanRefSol}}
%    \subfigure[Mean Approximate Solution ]{\includegraphics[width=0.47\textwidth]{approxsol}\label{MeanApproxSol}}
%  \end{center}
%
%  \caption{\label{meansol}{Comparison between mean reference solution of the
%fine-scale problem with the mean approximate solution obtained from the global-local multiscale approach using random values of the online parameter $\mu$.}}
%
%\end{figure}

\begin{figure}[htp]
  \begin{center}
    \includegraphics[width=0.47\textwidth]{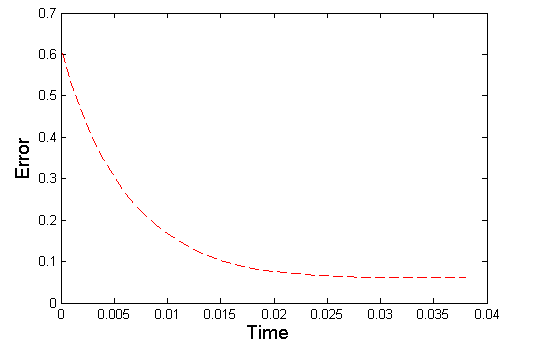}

  \end{center}

  \caption{\label{meanerror}{Mean error of approximating the solution by using global-local multiscale approach with random values of the online parameter $\mu$.}}

\end{figure}

\end{exmp}

\section{Conclusions}
\label{sec:conc}
In this work, we present a global-local nonlinear model reduction
approach to reduce the computational cost for solving high-contrast
nonlinear parabolic PDEs. This is achieved through two main stages;
offline and online. In the offline step, we use the generalized
multiscale finite element method (GMsFEM) to represent the coarse-grid
solutions through applying the local discrete empirical
interpolation method (DEIM) to approximate the nonlinear functions
that arise in the residual and Jacobian. Using the snapshots of the
coarse-grid solutions, we compute the proper orthogonal decomposition
(POD) modes. In the online step, we project the governing equation on
the space spanned by the POD modes and use the global DEIM to
approximate the nonlinear functions. Although one can perform global
model reduction independently of GMsFEM, the computations of the
global modes can be very expensive. Combining both local and global
mode reduction methods along with applying DEIM to inexpensively
compute the nonlinear function can allow a substantial speed-up.  We
demonstrate the effectiveness of the proposed global-local nonlinear
model reduction method on several examples of nonlinear multiscale
PDEs that are solved using a fully-implicit time marching schemes. The
results show the great potential of the proposed approach to reproduce
the flow field with good accuracy while reducing significantly the
size of the original problem.  Increasing
%The effect of the variations of
the number of the local and global modes  to improve the accuracy of the approximate
solution is examined. Furthermore, the robustness of proposed model reduction approach with respect to
variations in initial conditions, permeability fields, nonlinear-function's parameters, and forcing terms is demonstrated.

\end{document}